\documentclass{article}

\usepackage{amsmath}
\usepackage[dvips]{graphicx}
\usepackage{amsfonts}
\usepackage{amsopn}
\usepackage{amsthm}
\usepackage{amssymb}
\usepackage{url} 

\newtheorem{theorem}{\bf \normalsize \bf  Theorem}[section]
\newtheorem{definition}[theorem]{\bf \normalsize \bf Definition}

\newtheorem{lemma}[theorem]{\bf \normalsize \bf Lemma}

\renewcommand{\O}{\Omega}

\newcommand{\cO}{$O~$}

\newcommand{\G}{\Gamma}
\newcommand{\gr}{\nabla^{\prime}}

\renewcommand{\d}{\delta}

\newcommand{\pr}{\partial}

\newcommand\Sn{\mathbb {S}^{n}}

\newcommand\Rnn{\mathbb {R}^{n}}
\newcommand\Rn{\mathbb {R}^{n+1}}

\newcommand{\p}{\rho}

\newcommand{\n}{\nabla}

\newcommand{\np}{\n^{\prime}}

\newcommand{\Ob}{\bar{\Omega}}

\newcommand{\pb}{\bar{\psi}}

\newcommand{\la}{\lambda}

\newcommand{\gh}{\hat{g}}
\newcommand{\bh}{\hat{b}}
\newcommand{\ah}{\hat{a}}
\newcommand{\ub}{\bar{u}}
\newcommand{\vt}{\tilde{v}}
\def\proof{\noindent {\bf Proof.~}}
\def\endproof{\ \hfill QED. \bigskip}

\begin{document}

\title{Uniqueness of Starshaped Compact Hypersurfaces With Prescribed $m$-th 
Mean Curvature in Hyperbolic Space \footnotetext{1991 Mathematics Subject Classification 53A10, 35J60}}
\author{Jo\~ao Lucas M. Barbosa \footnotemark[1] \and Jorge H. S. de Lira \thanks {Partially supported by CNPq and FUNCAP} \\ \and \\ Vladimir Oliker \thanks{Partially supported by the National Science Foundation grant DMS-04-05622}}

\date{}

\maketitle

\begin{abstract}
Let $\psi$ be a given function defined on a Riemannian space.
Under what conditions does there exist a compact starshaped hypersurface $M$ for which $\psi$, when evaluated on $M$, coincides with the $m-$th elementary symmetric function of principal curvatures of $M$ for a given $m$? The corresponding existence and uniqueness problems in Euclidean space have been investigated by several authors in the mid 1980's. Recently, conditions for existence were established in elliptic space and, most recently, for hyperbolic space. However, the uniqueness problem has remained open. In this paper we investigate the problem of uniqueness in hyperbolic space and show that uniqueness (up to a geometrically trivial transformation) holds under the same conditions under which existence was established.
%Consider a positive real function defined in an annulus centered at certain point of the hyperbolic space. Assume %this function depends only on the distance from that point, satisfying suitable conditions on the boundary of the %annulus and whose quotient with curvature of the spheres centered at point is decreasing. Consider the problem of %finding the starshaped hypersurfaces with respect to the point, contained into the annulus whose m-th curvature is %equal to the given function. Recent it has been proved the existence of such hypersurfaces with the contribution of %several authors. We show that the solution is unique.    
\end{abstract}

\section{Introduction} \label{intro}
In Euclidean space $\Rn$ fix a point \cO and let $\Sn$ be the unit
sphere  centered at \cO. Let $u$ denote a point on
$\Sn$ and let $(u,\p)$ be the spherical coordinates in $\Rn$ with the origin at $O$. The
standard metric on $\Sn$ induced from $\Rn$ we denote by $e$.
Let $I = [0,a)$, where $a=const, ~0 < a  \leq \infty,$ and $f(\p)$ a positive
$C^{\infty}$ function on $I$ such that $f(0) =0$. Introduce in
 $\Sn \times I$ the metric 
\begin{equation}\label{metric}
h= d\p^2 + f(\p)e
\end{equation}
and consider the resulting Riemannian space.  When $a = \infty$ and $f(\p) = \p^2$
this space is the Euclidean space $\Rn\equiv \Rn(0)$, when $a = \infty$ and 
$f(\p) = \sinh ^2 \p$ it is the hyperbolic space $\Rn(-1)$ with sectional curvature $-1$ 
and when $a = \pi/2,~f(\p) = \sin ^2 \p$, it is the elliptic  space $\Rn(1)$
%$\Sn_+$ 
with sectional curvature $+1$. We use the notation $\Rn(K), ~K=0,\pm 1$ for either
of these spaces. 

Let $M$ be a hypersurface in $\Rn(K)$ and $m$, $1\leq m\leq n$,
be an integer.
The $m$-th mean curvature,  $H_m(\la)\equiv H_m(\la_1,...,\la_n)$, of $M$ is the normalized
elementary symmetric function of order $m$ of the principal curvatures 
$\la_1,...,\la_n$ of $M$, that is,
\[H_m(\la) = (^n_m)^{-1}\sum_{i_1<...<i_n}\la_{i_1}\cdots \la_{i_m}.\]

The subject of this paper is the following problem.
Let $\psi(u,\p),~~u \in \Sn,~~\p \in I,$ be a given positive function and
$m, ~1\leq m\leq n,$ a given integer. Under what
conditions on $\psi$ does there exist a smooth hypersurface $M$ in 
$\Rn(K)$ given
as $(u,z(u)), ~~u \in \Sn,~~z > 0,$ for which 
\begin{equation} \label{eq_0}
H_m (\la(z(u))) = \psi (u,z(u))~\forall u \in\Sn?
\end{equation}
In addition, if such a hypersurface
exists then we wish to know conditions for uniqueness.

In analytic form this problem consists in establishing existence and uniqueness
of solutions for a second order nonlinear partial differential equation on $\Sn$ expressing $H_m$ in terms of $z$. When $m=1$ this equation is quasilinear and for $m >1$ it is fully nonlinear. In particular, when $m=n$ it is of Monge-Amp\`{e}re type. In Euclidean space $\Rn(0)$
this problem was investigated and conditions for existence and uniqueness were
given 
by I. Bakelman and B. Kantor \cite{BK1, BK2}
and A. Treibergs and S.W. Wei \cite{TW} when $m = 1$
(the mean curvature case), by V. Oliker \cite{O:curv} when $m = n$
(the Gauss curvature case), and by L. Caffarelli, L. Nirenberg and
J. Spruck \cite{CNS} when $1 < m < n$.
 
In  \cite{O:curv:89} V. Oliker
investigated  the problem for hypersurfaces in $\Rn(-1)$
and $\Rn(1)$ when $m=n$ and gave conditions for existence and uniqueness.  
In \cite{BLO} L. Barbosa, J. Lira and V. Oliker
obtained $C^0,~~C^1$ and $C^2$
estimates for solutions 
of (\ref{eq_0})  
for the elliptic space form $\Rn(1)$ for any $m,~1\leq m \leq n,$ and then, in \cite{LO}, Y. Y. Li and V. Oliker,
used these estimates and degree theory for fully nonlinear elliptic operators 
\cite{YYLi_deg} to prove existence of solutions. In the same paper \cite{BLO}, the authors also obtained the $C^0$ and $C^1$ estimates for any 
$m,~1\leq m \leq n,$ in
the hyperbolic space $\Rn(-1)$. Recently, Q. Jin and Y. Y. Li \cite{JL} obtained the $C^2$ estimates for $\Rn(-1)$ and
proved existence for this case as well. The main results in \cite{LO} and
\cite{JL} can be formulated together as follows. 

Denote by $\G_m$ the connected component
of $\{\la \in \Rnn~|~H_m(\la)> 0\}$ containing the positive cone
$\{\la \in \Rnn~|~\la_1,...,\la_n > 0\}.$
\begin{definition} A positive function $z\in C^2(\Sn)$ is $m-$admissible for
the operator $H_m$  if the corresponding hypersurface
$M=(u,z(u)),~~u \in \Sn$, is such that
at every point of $M$ 
the principal curvatures $(\la_1(z(u)),...,\la_n(z(u))) \in \G_m$,
where the $\la_i$ are calculated with respect to the inner normal.
\end{definition}
\begin{theorem}\label{MT}  Let $1\le m \le n, ~K =\pm1,$ and
$\psi (u,\p)$  is a positive smooth function
on the annulus $\bar{\O} \subset \Rn(K)$, $\bar{\O}: ~u \in \Sn,~~\p
\in [R_1,R_2]$, where
$0 < R_1 <  R_2 < a$, and $a=\infty $ for $\Rn(-1)$ and 
$a= \pi/2$ for $\Rn(1)$. 
Suppose $\psi$ satisfies the conditions:

If $K=-1$
\begin{eqnarray}
\psi(u,R_1) \geq \coth ^m R_1 ~
\mbox{for}~u \in \Sn, \label{cond1} \\
\psi(u,R_2) \leq \coth ^m R_2~
\mbox{for}~~u \in \Sn, \label{cond2}
\end{eqnarray}
and
\begin{equation} \label{cond3}
\frac{\pr}{\pr \p}\left [\psi(u,\p)\sinh ^{m}\p \right]\le 0
~~\mbox{for all}~u \in \Sn~~\mbox{and}~
\p \in [R_1,R_2];
\end{equation}
\indent If $K=1$
\begin{eqnarray}
\psi(u,R_1) \geq \cot ^m R_1 ~
\mbox{for}~u \in \Sn, \label{cond1_2} \\
\psi(u,R_2) \leq \cot ^m R_2~
\mbox{for}~~u \in \Sn, \label{cond2_2}
\end{eqnarray}
and
\begin{equation} \label{cond3_2}
\frac{\pr}{\pr \p}\left [\psi(u,\p)\cot ^{-m}\p \right]\le 0
~~\mbox{for all}~u \in \Sn~~\mbox{and}~
\p \in [R_1,R_2].
\end{equation}
Then there exists a closed, smooth,
embedded hypersurface $M$ in $\Rn(K)$, $M \subset \Ob$, which
is a radial graph over $\Sn$ of an $m-$admissible function $z$
and
\begin{equation} \label{eq0}
H_m (\la(z(u))) = \psi (u,z(u)) ~~~\mbox{for all}~~  u \in \Sn.
\end{equation}
\end{theorem}

Similar to the case of $\Rn(1)$  the proof in \cite{JL} uses degree theory. The degree theory arguments
in \cite{LO} and \cite{JL} do not provide an answer to the uniqueness problem
and thus for the elliptic and hyperbolic space forms this question remained open
except for the case $m=n$ \cite{O:curv:89}.
The purpose of this paper is to show that under the same conditions as in
Theorem \ref{MT} we can also
prove uniqueness for the hyperbolic space for all $m,~1\leq m \leq n$. 
Namely, we have the following
\begin{theorem} \label{MT1}
Let $K=-1$. Then under condition  (\ref{cond3})  in Theorem \ref{MT}
any two hypersurfaces defined by $m-$admissible solutions $z^1$ and $z^2$ 
of (\ref{eq0}) in
$\Ob$  are related by the transformation:
\begin{equation}\label{K-scaling1}
c\tanh(\frac{z_1(u)}{2})= \tanh(\frac{z_2(u)}{2}),~~u \in \Sn,
\end{equation}
where $c$ is a positive constant.
If the inequality (\ref{cond3})  is strict then $c=1$, that is,
the hypersurface $M$ in Theorem \ref{MT} is unique. 
\end{theorem}

For $m=n$ the condition (\ref{cond3}) is slightly less restrictive than condition c) in Theorem 1.1 in \cite{O:curv:89}.
%In a certain sense, this completes the solution of the above stated problem for 
%space forms.
For 
the elliptic space the uniqueness problem is still open except for the already known case when $m=n$. In this case condition (\ref{cond3_2}) also implies uniqueness.

\section{\protect The Equation of the Problem }\label{loc}
In this section we present some local formulas and lemmas valid in $\Rn(K)$
where $K=\pm 1$. Though our main result (Theorem \ref{MT}) applies only to the
case $K=-1$, it seems worthwhile to record here the results which are also valid for the case $K=+1$ because they should
 be useful in future studies of similar problems. Furthermore, the presentation in this section is carried out  in a unified way simultaneously for both cases.  

{\bf 1. The main equation.} First we fix our notation.
Unless explicitly stated otherwise, the range for the latin
indices is $1,...,n$. The summation convention over repeated lower and upper indices
is assumed to be in effect. 
Denote by  $(u^1,...,u^n)=u$ smooth local coordinates on $\Sn$ and let
$\pr_i = \pr/\pr u^i, ~i=1,2,...,n,$ be the corresponding
 local  frame of tangent vectors  such that $e(\pr_i,\pr_j)=e_{ij}$. 
The first covariant derivative of a function $v\in C^2(\Sn)$ is given by
$v_i \equiv \gr_{i}v = \pr_i v$. Put
$(e^{ij}) = (e_{ij})^{-1}$ and let
$$\gr v = v^i\pr_i, ~~\mbox{where}~~v^i = e^{ij}v_j.$$
For the covariant derivative of $\gr v$ we
have
\[\gr_{\pr_s}\gr v = v_{sj}e^{ji}\pr_i + v_j\gr_{\pr_s}(e^{ij}\pr_i) =
(v_{sj}  - \G_{sj}^{\prime i}v_i) e^{jk}\pr_k,\]
where
\[
v_{sj} = 
\frac {\pr^2 v}{\pr u^s \pr u^j}\]
and $\G_{sj}^{\prime i}$ are the Christoffel symbols of the second kind
of the metric $e$. 
The second covariant derivatives of $v$ are defined by
\begin{equation}\label{covdev2}
\gr_{sj}v = v_{sj}  - \G_{sj}^{\prime i}v_i.
\end{equation}

Next we recall some of the basic formulas
derived in \cite{BLO}. 
Let $M$
be a hypersurface in $\Rn (K)$ given by $r(u) = (u,z(u)), ~u \in \Sn$, where $z \in C^2(\Sn)$ 
and 
positive on $\Sn$. 
The metric  $g = g_{ij}du^idu^j$  induced on $M$ from $\Rn(K)$ 
has coefficients
\begin{equation}\label{form1:prel}
g_{ij} = fe_{ij} + z_iz_j~~\mbox{ and}~~ 
\det(g_{ij})= f^{n-1}(f+|\gr z|^2)\det(e_{ij}).
\end{equation}

The elements of the inverse matrix $(g^{ij}) = (g_{ij})^{-1}$ are
\begin{equation}\label{ginv}
g^{ij} = \frac{1}{f}\left [ e^{ij} - \frac{z^iz^j}{f+|\gr z|^2}\right ].
\end{equation}

With the choice of the normal on $M$ in inward direction the second fundamental form $b$ of $M$ has coefficients:
\begin{equation}\label{2form}
b_{ij} = \frac{f}{\sqrt{f^2 + f|\gr z|^2}}\left[-\gr_{ij}z + 
\frac {\pr \ln f}{\pr \p}z_iz_j + 
\frac{1}{2}\frac {\pr f}{\pr \p}e_{ij}\right].
\end{equation}
Note that the second fundamental form
of a sphere $z = const > 0$ is positive definite, since for  $\Rn(K)$
$\pr f/\pr \p > 0$. 

The principal curvatures of $M$ at a point $(u,z(u))$ are the eigenvalues of the second
fundamental form $b$ relative to the metric $g$  and are  the real roots,
$\la_1(z(u)),...,\la_n(z(u))$, of the
equation
\[\det(b_{ij}(z(u))-\la g_{ij}(z(u))=0\]
or, equivalently, of
\[\det(a_j^i(z(u)) - \la \d_j^i) = 0,\]
where 
\begin{equation}\label{a-def}
a_j^i = g^{ik}b_{kj},
\end{equation}
is a self-adjoint transformation of the tangent space to $M$ at $(u,z(u))$. 
The elementary symmetric function of order $m,~ 1 \le m \le n$, 
of the principal curvatures
%$\la(z(u)) =(\la_1(z(u),...,\la_n(z(u))$
is
\begin{equation}\label{sum_of_minors}
S_m(\la) = \sum_{i_1<...<i_n}\la_{i_1}\cdots \la_{i_m}~~~\mbox{and}~~~
S_m(\la) = (^n_m)H_m(\la) = F_m(a_j^i),
\end{equation}
where $F_m$ is the sum of  principal minors of $(a_j^i)$ of
order $m$.  
Evidently,
\begin{equation}\label{eqF}
F_m(a_j^i(z(u))) \equiv F(u,z, \gr_1,...,\gr_n z, \gr_{11}z,...,\gr_{nn}z),
\end{equation}
and the equation (\ref{eq0}) assumes the form
\begin{equation}\label{eq1}
S_m(\la(z(u)))\equiv F_m(a_j^i(z(u))) = \pb (u,z(u)),
\end{equation}
 where 
 $\pb \equiv (^n _m)\psi$.

{\bf 2. The conformal model of $\Rn(K)$ and a change of the function $z$.} 
For the function $f(\p), ~\p\in I,$ in (\ref{metric}) corresponding to
$\Rn(-1)$ or $\Rn(+1)$ we put $$s(\p)=\sqrt{f(\p)},~~c(\p) = \frac{d s(\p)}{d \p}, ~~t(\p)=\frac{s(\p)}{c(\p)}.$$

It will be convenient to make a change of the function $z$ in (\ref{eqF})
by setting
$v(u) = t(z(u)/2)$\footnote{ This is equivalent to re-writing (\ref{metric}) in the conformal model of the corresponding space form in the unit ball in Euclidean space $\Rn$ centered at the origin.}. Put
\[q = \frac{2}{1+Kv^2}.\]
Then
\begin{equation}\label{z-to-v}
z_i = qv_i, ~~~\gr_{ij}z= q\gr_{ij}v-Kq^2vv_iv_j.
\end{equation}
Put
\[\gh_{ij}(v)=v^2e_{ij}+v_iv_j,~\gh^{ij}(v)=\frac{1}{v^2} \left ( e^{ij} - \frac{v^iv^j}{W^2(v)} \right), ~W(v) =  \sqrt{v^2 +|\gr v|^2}.\]
A substitution into (\ref{ginv}) gives
\[g^{ij}(v)=\frac{1}{q^2}\gh^{ij}(v)\]
and a substitution into (\ref{2form}) gives
\[b_{ij}(v)=q\bh_{ij}(v)-Kq^2v^2\frac{\gh_{ij}(v)}{W(v)},\]
where
\begin{equation}\label{2-nd form}
\bh_{ij}(v)=\frac{-v\gr_{ij}v+2v_iv_j+v^2e_{ij}}{W(v)}.
\end{equation}
Note that $\gh$ and $\bh$ are respectively the first and second fundamental
forms in the Euclidean sense of the hypersurface which is a graph of $v$ over $\Sn$ in the unit ball \cite{O:curv}.
Finally, we obtain
\begin{equation} \label{curv-1}
a^i_j(v)=g^{ik}(v)b_{kj}(v)= \frac{\ah^i_j(v)}{q}-K\frac{v^2\d^i_j}{W(v)},
~~\mbox{where}~~\ah^i_j(v)=\gh^{ik}(v)\bh_{kj}(v).
\end{equation}
For a $m-$admissible function $z \in C^2(\Sn)$
and $v=t(z/2)$ consider the family of functions  $sv$, where $s > 0$ and such that
$sv < 1$. Then
\begin{equation}\label{curv-s}
a^i_j(sv)= \frac{1+Ks^2v^2}{2s}\ah^i_j(v)-K\frac{sv^2}{W(v)}\d^i_j.
\end{equation}

Define, as before, the eigenvalues 
$\lambda_i(sv(u)),i=1,...,n,$ of  $(b_{ij}(sv(u)))$ with
respect to $(g_{ij}(sv(u)))$ (which is positive definite) and consider
the corresponding $m-$th elementary symmetric function $S_m(\la(sv(u)))$.  Clearly, since $z$ is $m-$admissible,
the function $v$ is $m-$admissible, that is,  $\la(v(u))\in \G_m$.
\begin{lemma}\label{sv}
Let $z, v$ and $s$ be as above. 
Put
\[A(sv)=\frac{1+Ks^2v^2}{s(1+Kv^2)},~B(sv)=K\frac{(1-s^2)v^2}{s(1+Kv^2)W(v)}.\]
 Then
\begin{equation}\label{sv-0}
S_m(\lambda(sv))=
A^m(sv)S_m(\lambda(v))+\sum_{j<m}c(n,m,j)A^j(sv)B^{m-j}(sv)S_j(\la(v)),
\end{equation}
where $c(n,m,j)$ are positive coefficients.
Furthermore,
if $K=-1$ and $s \geq 1$ or if $K=+1$ and $s \leq 1$, then
\begin{equation} \label{sv-1}
S_m(\lambda(sv(u)))\geq
A^m(sv(u))S_m(\lambda(v(u))).
\end{equation}
In particular, $sv$ is $m-$admissible for $H_m$ in $\Rn(K)$ for the corresponding choice of $s$. 
\end{lemma}
\proof 
It follows from (\ref{curv-1}) and (\ref{curv-s}) that
\[
a^i_j(sv)=A(sv)a^i_j(v)+B(sv)\d^i_j.
\]
Since $v$ is $m-$admissible, $S_j(v) > 0$ for each $j \leq m$ (see \cite{CNS3})
and $A(sv) > 0$ because $sv < 1$. On the other hand, $B(sv) \geq 0$ with each 
 choice
of $s$ as in the statement of the lemma. Then $S_m(\lambda(sv)) >0$ in both cases.
Because $sv$ is a positive multiple of $v\in \G_m$ we conclude that $sv \in \G_m$
in both cases.
\endproof

We complete this section with the following
\begin{lemma} \label{negell}Let $z$ be $m-$admissible for the operator $H_m$ and $v =t(z/2)$. Then the operator $F_m(a^i_j(v))$ is negatively elliptic on $v$ on $\Sn$.
\end{lemma}
\begin{proof} In order to show that $F_m(a^i_j(v))$ is negatively elliptic 
we need to show that at any point of $\Sn$ the quadratic form
\begin{equation}\label{quadratic}
\frac{\pr F_m(a_j^i(v))}{\pr \gr_{ij}v}\xi^i\xi^j < 0 ~~\forall \xi \in \Rnn, 
~\xi \neq 0.
\end{equation}
It folows from (\ref{curv-1}) and (\ref{2-nd form}) that
\begin{equation}\label{quadratic1}
\frac{\pr F_m}{\pr \gr_{ij}v} = -\frac{v}{q(v)W(v)}F_i^j,~~\mbox{where}~~ F_i^j=\frac{\pr F_m}{\pr a^i_j(v)}.
\end{equation}
Thus, we need to consider the matrix $(F_i^j)$. 

Fix an arbitrary point $u_0 \in \Sn$ and diagonalize at that point the
metric $(g_{ij}(v))$ and the second fundamental form $(b_{ij}(v))$ using the orthonormal set of principal directions as a basis. Then at $u_0$ we have
$g_{ij}(v) = \delta_{ij}$,
\[
b_{ij}(v) = \left \{
\begin{array}{ll}
\la_i(v) &\mbox{when}~~i=j\\
0 &\mbox{when} ~~i\neq j,
\end{array}
\right.
\]
and
\[
a^i_j(v) = \left \{
\begin{array}{ll}
\la_i(v) &\mbox{when}~~i=j\\
0 &\mbox{when} ~~i\neq j.
\end{array}
\right.
\]
Then at $u_0$ we have
\begin{equation} \label{quadratic2}
F^j_i = 0 ~\mbox{when}~i\neq j~\mbox{and}~F^i_i=\frac{\pr S_m(\la_1(v),...,\la_n(v))}{\pr \la_i(v)},
\end{equation}
where no summation over $i$ is performed. Since $z$, and therefore $v$, are
$m-$admissible for $H_m$, it follows that
$S_m(\la_1(v),...,\la_n(v)) > 0$. Then, by a well known property of
elementary symmetric functions $\frac{\pr S_m}{\pr \la_i(v)} > 0$ 
for each $i=1,2,...,n$; see \cite{CNS3}.
Now, (\ref{quadratic}) follows from (\ref{quadratic2}) and (\ref{quadratic1}).
\end{proof}

\section{\protect{ Proof of  Theorem \ref{MT1}}}\label{mt}
In this section we work in the hyperbolic
space $\Rn(-1)$.

Let $z_1$ and $z_2$ be two different $m-$admissible solutions of (\ref{eq0}) and $M_1, M_2$ the corresponding hypersurfaces on the annulus
$\Ob$. It follows from Lemma \ref{lemmaA} (see appendix)
that for any
 $m-$admissible solution $z$ of (\ref{eq0}) such that
$R_1\leq z(u)\leq R_2$
we have either $z(u) \equiv R_1$ or  $z(u) \equiv R_2$ or 
\begin{equation} \label{assumpt1}
R_1 < z(u) < R_2 ~~\forall u \in \Sn.
\end{equation}
Assume first that 
\begin{equation} \label{hyp1}
R_1 \leq z_k(u) < R_2 ~~\mbox{ for} ~k=1,2~\mbox{and}~\forall u \in \Sn.
\end{equation}
The case  when  $z_1$ or $z_2 \equiv R_2$ is special and
will be   treated separately at the end of the proof.

%It will be convenient to divide the proof into two cases: $K=-1$ and $K=+1$. {\bf %We begin with the  case $K=-1$}.

Let $v_k(u) = t(z_k(u)/2)$, where now $t(z_k(u)/2)=\tanh(z_k(u)/2)$. 
Suppose $v_1 < v_2$ somewhere on $\Sn$; otherwise re-label them.  Multiply $v_1$ by $s \geq 1$ such that
 $$sv_1(u) < 1, ~sv_1(u)\geq v_2(u) ~\forall u \in \Sn~~\mbox{and}~~
  sv_1(\ub)=v_2(\ub)~~\mbox{at some}~~\ub \in \Sn.$$ By
(\ref{hyp1}) there exists some neighborhood $U\subset \Sn$ of the point $\ub$
such that $z^s=2t^{-1}(sv_1)$ satisfies the inequality
\[R_1 < z^s(u) < R_2,~u\in U.\]
Since $S_m(\lambda(z_1))=\pb(u,z_1(u))$, it follows from Lemma \ref{sv} that in $U$
\begin{equation}\label{monot}
S_m(\lambda(sv_1))-\pb(u,2t^{-1}(sv_1)) \geq A^m(sv_1)\pb(u,2t^{-1}(v_1))-\pb(u,2t^{-1}(sv_1)).
\end{equation}
Put $sv_1=\tilde{v}$. 
Then, using the explicit expression for $A(sv_1)$ and taking into account that $K=-1$, we get
\begin{equation}\label{monot1}
S_m(\lambda(\vt))-\pb(u,2t^{-1}(\vt)) \geq \left[\frac{1-\vt^2}{s(1-\frac{\vt^2}{s^2})}\right ]^m\pb(u,2t^{-1}(\frac{\vt}{s}))-\pb(u,2t^{-1}(\vt)).
\end{equation}
Put
\[Q(s)=\left[\frac{1-\vt^2}{s(1-\frac{\vt^2}{s^2})}\right ]^m\pb(u,2t^{-1}(\frac{\vt}{s}))-\pb(u,2t^{-1}(\vt)).\]
Note that $Q(1)\equiv 0$. We have
%\begin{eqnarray}
$$\frac{\partial Q}{\partial s} =  
 -m\left[\frac{1-\vt^2}{s(1-\frac{\vt^2}{s^2})}\right ]^{m-1}
\frac{(1-\vt^2)(1+\frac{\vt^2}{s^2})}{s^2(1-\frac{\vt^2}{s^2})^2}\pb(u,2t^{-1}
(\frac{\vt}{s}))-$$ 
%\nonumber \\
$$\left[\frac{1-\vt^2}{s(1-\frac{\vt^2}{s^2})}\right ]^m\frac{2\vt}{s^2(1-\frac{\vt^2}{s^2})}\pb_z(u,2t^{-1}(\frac{\vt}{s}))=
$$
% \nonumber \\
$$-\frac{1}{s^{m+1}}\left [\frac{1-\vt^2}{1-\frac{\vt^2}{s^2}}\right ]^m
\left[m\frac{1+\frac{\vt^2}{s^2}}{1-\frac{\vt^2}{s^2}}\pb(u,2t^{-1}(\frac{\vt}{s}))
 +\frac{2\vt}{s(1-\frac{\vt^2}{s^2})}\pb_z(u,2t^{-1}(\frac{\vt}{s}))\right ]\geq 0,
$$
%\nonumber
%\end{eqnarray}
where $\pb_z=\frac{\partial \pb}{\partial z}$.
The last inequality on the right follows from  (\ref{cond3}). 

By (\ref{monot1}), (\ref{monot}) and the assumption
that $z_2$  is an $m-$admissible solution of (\ref{eq0}) we have
\[F_m(a^i_j(\vt))-\pb(u,2t^{-1}(\vt)) \geq 0=F_m(a^i_j(v_2))-\pb(u,2t^{-1}(v_2)).\]
Since $F_m$ is negatively elliptic, $\vt \geq v_2$ in $U$ and $\vt(\ub)=  v_2(\ub)$,
we conclude from the geometric form of
Aleksandrov's maximum principle \cite{A} that $\vt \equiv v_2$ in
$U$. By continuity, the set 
$$\{u\in \Sn~|~\vt(u) = v_2(u)\}$$ is open and closed on $\Sn$. Hence, $\vt(u)  = v_2(u)=sv_1(u)$ everywhere on $\Sn$ and the proof of uniqueness is complete in this case.

%In addition, the left hand side in (\ref{monot1}) vanishes and the function
%$Q(s)$ is nondecreasing and vanishes at $s=v_2/v_1$. Then $Q(s^{\prime})\equiv 0~
%\forall s^{\prime}\in [1,s]$.

Suppose now that $z_2 \equiv R_2$ and $z_1 < R_2~\forall u \in \Sn$. In this case we extend $\pb(u,\p)$ smoothly for $\p > R_2$ satisfying conditions 
\begin{eqnarray}
\psi(u,\p) \leq \coth ^m R_2~
\mbox{for}~~u \in \Sn, \label{cond2prime}
\end{eqnarray}
and
\begin{equation} \label{cond3prime}
\frac{\pr}{\pr \p}\left [\psi(u,\p)\sinh ^{m}\p \right]\le 0
~~\mbox{for all}~u \in \Sn~~\mbox{and}~
\p \geq R_2;
\end{equation}
Then, again, the same arguments apply and this completes the proof of the theorem.

\section{\protect{Appendix}}\label{lemmaA0}
\begin{lemma}\label{lemmaA} Assume that the conditions of the Theorem \ref{MT} are satisfied
except for conditions (\ref{cond3}) and (\ref{cond3_2}).  Then a
 $m-$admissible solution $z$ of (\ref{eq0}) such that
$R_1\leq z(u)\leq R_2$
is either $\equiv R_1$ or  $\equiv R_2$, or 
\begin{equation} \label{assumpt2}
R_1 < z(u) < R_2 ~~\forall u \in \Sn.
\end{equation}
\end{lemma}
This lemma was stated in \cite{BLO} without a detailed proof. At the suggestion
of the referee we provide a  proof here. The proof consists in
showing that the conditions of Aleksandrov's maximum principle \cite{A}
are satisfied.

\begin{proof}
Suppose, on the contrary, that there exists some $u_0 \in \Sn$ such that
$z(u_0)=R_2$ and $z(u)\not \equiv R_2$. (The case 
when $z(u_0)=R_1, z(u)\not \equiv R_1,$ is treated similarly.) Then $z$ attains a maximum at $u_0$. 
Consider the family of functions
\[
z(s) = (1-s)z + sR_2,~s\in [0,1].
\]
Obviously, $z(s)$ also attain a maximum $=R_2$ at $u_0$ for all $s\in [0,1]$.
We will need
an expression for the $m-$th elementary symmetric function
of the hypersurface $M(s)$ defined by $z(s)$ at $u_0$. We
have 
\begin{equation}\label{lemma_eq1}
\np z(s) = 0~\mbox{and}~\np_{ij} (z(s)) =(1-s) z_{ij} ~\mbox{at}~u_0,~~s\in [0,1].
\end{equation}
Put
\[
\mu= \frac{1}{2f(R_2)}\frac{\pr f(z(s))}{\pr z}|_{z(s)=R_2}
\]
and observe that $\mu > 0$, since $0 < R_2 < a$ and $\frac{\pr f}{\pr \p} > 0$
in $\Ob$.
Using (\ref{ginv}), (\ref{2form}) and (\ref{covdev2}) and noting that
\begin{equation*}
a^i_j(R_2)= \mu\d^i_j,
\end{equation*}
we obtain at the point $u_0$
\begin{eqnarray*}
a^i_j(z(s)) = 
\frac{e^{ik}}{f(R_2)}\left[-(1-s)z_{kj}+\mu e_{kj} \right]=
(1-s)a^i_j(z)+sa^i_j(R_2).
\end{eqnarray*}
Then
\begin{equation}\label{lemma_eq2}
S_m(\la(z(s)))|_{u_0}= \sum_{p=0}^m (1-s)^p(\mu s)^{m-p}S_p(\la(z))|_{u_0},
\end{equation}
where $S_0=(^n_m)$. Since $S_p(\la(z)) > 0$ for all $p \leq m$ and $\mu > 0$,
it follows that $S_m(\la(z(s)))|_{u_0} > 0$
for all $s \in [0,1]$. By continuity $S_m(\la(z(s))) > 0$ in some neighborhood
$U_0$ of $u_0$ in $\Sn$. Then by \cite{CNS3} $\frac{\pr S_m}{\pr \la_i(z(s))} > 0$
and by shrinking $U_0$, if necessary, we have 
$\frac{\pr S_m}{\pr \la_i(z(s))} \geq C > 0$ for all $u \in U_0$ with some
fixed constant $C$. It follows now from Lemma  \ref{negell} and the second
expression in (\ref{z-to-v}) that $-F_m(a^i_j(z(s))$
is positively elliptic in $U_0$ for all $s \in [0,1]$. 

The above arguments establish that the function 
$-F_m(a^i_j(z(s))+\pb(u,z(s))$
satisfies the conditions (1)-(4) in \cite{A}, \S1. 
(Note that our orientation of $M(s)$ is opposite to that in \cite{A}.) We need to check one more
inequality. Namely, since $z$ satisfies (\ref{eq0}) on $\Sn$ we have
\[-F_m(a^i_j(z)+\pb(u,z)=0,\]
while, taking into account (\ref{cond2}) for the hyperbolic space or
(\ref{cond2_2}) for the elliptic space, we also have
\[-F_m(a^i_j(R_2)+\pb(u,R_2)) \leq 0.\]
Since, also, $z(u) \leq R_2$ on $\Sn$ and $z(u_0) = R_2$ it follows from
the maximum principle in \cite{A} that $z(u) = R_2$ everywhere on $U_0$.
This implies that the set $\{u\in \Sn~|~z(u)=R_2\}$ is open on $\Sn$.
Since it is also closed, we conclude that $z(u) = R_2$ everywhere on $\Sn$. 
\end{proof}

\noindent{\bf Address} \\
\noindent Jo\~{a}o Lucas M. Barbosa and Jorge  H.S. de Lira \\{\sc Universidade Federal do Cear\'{a}, Fortaleza, Brazil } \\
\noindent Vladimir I. Oliker \\ {\sc Emory University 
Atlanta, Georgia  USA}

\end{document}